\documentclass{egpubl}
\usepackage{sigrad2016}

 \WsPaper         
 \electronicVersion 

\ifpdf \usepackage[pdftex]{graphicx} \pdfcompresslevel=9
\else \usepackage[dvips]{graphicx} \fi

\PrintedOrElectronic

\usepackage{t1enc,dfadobe}

\usepackage{egweblnk}
\usepackage{cite}

\usepackage[capitalise,nameinlink]{cleveref}
\crefname{section}{Section}{Sections}
\crefname{figure}{Figure}{Figures}
\crefname{equation}{Equation}{Equations}

\usepackage[caption=false]{subfig}
\usepackage{bm}
\usepackage{mathtools}
\usepackage{amsmath}
\usepackage{xfrac}
\usepackage{multirow}
\usepackage{textcomp}

\title[RBF Approximation for Large Datasets]%
      {A Radial Basis Function Approximation for Large Datasets}

\author[Z. Majdisova \& V. Skala]
       {Z. Majdisova$^{1}$
        and V. Skala$^{1}$
        \\
         $^1$Department of Computer Science and Engineering, Faculty of Applied Sciences, University of West Bohemia,\\ Univerzitni 8, CZ 30614 Plzen, Czech  Republic
       }

\begin{document}


\maketitle

\begin{abstract}
Approximation of scattered data is often a task in many engineering problems. The Radial Basis Function (RBF) approximation is appropriate for large scattered datasets in $d$-dimensional space. It is non-separable approximation, as it is based on a distance between two points. This method leads to a solution of overdetermined linear system of equations.

In this paper a new approach to the RBF approximation of large datasets is introduced and experimental results for different real datasets and different RBFs are presented with respect to the accuracy of computation. The proposed approach uses symmetry of matrix and partitioning matrix into blocks.


\end{abstract}

\section{Introduction}
Interpolation and approximation are the most frequent operations used in computational techniques. Several techniques have been developed for data interpolation or approximation, but they mostly expect an ordered dataset, e.g. rectangular mesh, structured mesh, unstructured mesh etc. However, in many engineering problems, data are not ordered and they are scattered in $d-$dimensional space, in general. Usually, in technical applications the conversion of a scattered dataset to a semi-regular grid is performed using some tessellation techniques. However, this approach is quite prohibitive for the case of $d-$dimensional data due to the computational cost.

Interesting techniques are based on the Radial Basis Function (RBF) method which was originally introduced by \cite{Hardy1971}. They are widely used across of many fields solving technical and non-technical problems. The RBF applications can be found in neural networks, data visualization \cite{pepper2014}, surface reconstruction \cite{carr2001}, \cite{Turk2002}, \cite{pan2011two}, \cite{nedved2013}, \cite{skala2014}, solving partial differential equations \cite{Li2013}, \cite{Hon2015}, etc. The RBF techniques are really meshless and are based on collocation in a set of scattered nodes. These methods are independent with respect to the dimension of the space. The computational cost of these techniques increase nonlinearly with the number of points in the given dataset and linearly with the dimensionality of data.

There are two main groups of basis functions: global RBFs and Compactly Supported RBFs (CS-RBFs) \cite{wendland2006}. Fitting scattered data with CS-RBFs leads to a simpler and faster computation, but techniques using CS-RBFs are sensitive to the density of scattered data. Global RBFs lead to a linear system of equations with a dense matrix and their usage is based on sophisticated techniques such as the fast multipole method \cite{Darve2000195}. Global RBFs are useful in repairing incomplete datasets and they are insensitive to the density of scattered data.

For the processing of scattered data we can use the RBF interpolation or the RBF approximation. The RBF interpolation, e.g. presented by \cite{Skala2015}, is based on a solution of a linear system of equations:
\begin{equation}
\mathbf{Ac=h}\textrm{,}
\end{equation}
\noindent where $\mathbf{A}$ is a matrix of this system, $\mathbf{c}$ is a column vector of variables and $\mathbf{h}$ is a column vector containing the right sides of equations. In this case, $\mathbf{A}$ is an $N\times N$ matrix, where $N$ is the number of points in the given scattered dataset, the variables are weights for basis functions and the right sides of equations are values in the given points. The disadvantage of RBF interpolation is the large and usually ill-conditioned matrix of the linear system of equations. Moreover, in the case of an oversampled dataset or intended reduction, we want to reduce the given problem, i.e. reduce the number of weights and used basis functions, and preserve good precision of the approximated solution. The approach which includes the reduction is called the RBF approximation. In the following section, the method recently introduced in \cite{skala2013} is described in detail. This approach requires less memory and offer higher speed of computation than the method using Lagrange multipliers \cite{fasshauer2007}.
Further, a new approach to RBF approximation of large datasets is presented in the \cref{sec:big}. These approach uses symmetry of matrix and partitioning matrix into blocks.

\section{RBF Approximation}\label{sec:RBF}
\begin{figure}[htb]
\centering\includegraphics[width=.9\linewidth]{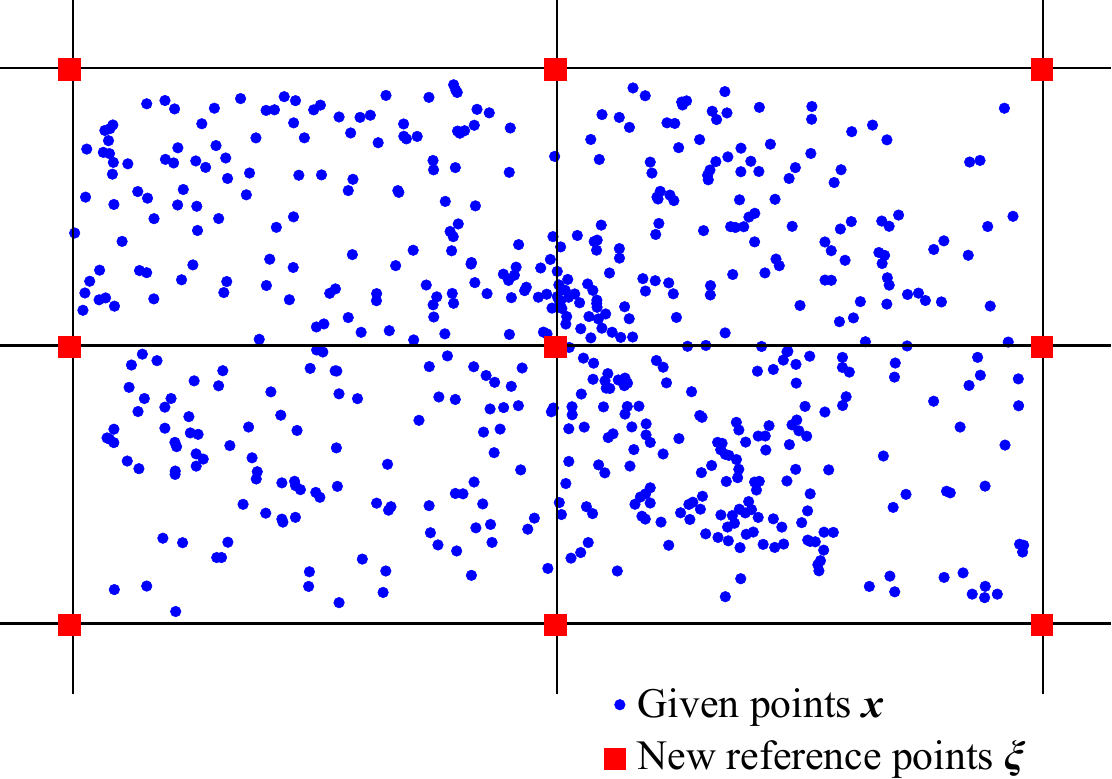}
\caption{The RBF approximation and reduction of points.}
\label{fig:virtual_point}
\end{figure}

For simplicity, we assume that we have an unordered dataset $\{\mathbf{x}_i\}_1^N \in E^2$. However, this approach is generally applicable for $d$-dimensional space. Further, each point $\mathbf{x}_i$ from the dataset is associated with a vector $\mathbf{h}_i\in E^p$ of the given values, where $p$ is the dimension of the vector, or scalar value, i.e. $h_i\in E^1$. For an explanation of the RBF approximation, let us consider the case when each point $\mathbf{x}_i$ is associated with a scalar value $h_i$,  e.g. a $2 \sfrac{1}{2}D$ surface. Let us introduce a set of new reference points $\{\bm{\xi}_j\}_1^M$, see \cref{fig:virtual_point}.

These reference points may not necessarily be in a uniform grid. It is appropriate that their placement reflects the given surface (e.g. the terrain profile, etc.) as well as possible. The number of reference points $\bm{\xi}_j$ is $M$, where $M \ll N$. Now, the RBF approximation is based on the distance computation of the given point $\mathbf{x}_i $ and the reference point $\bm{\xi}_j$.

The approximated value is determined similarly as for interpolation (see \cite{Skala2015}):
\begin{equation}
\label{eq:approximant}
\displaystyle f(\mathbf{x})=\sum_{j=1}^M c_j \phi(r_j)=\sum_{j=1}^M c_j \phi(\|\mathbf{x}-\bm{\xi}_j\|)\textrm{,}
\end{equation}
\noindent where $\phi(r_j)$ is a used RBF centered at point $\bm{\xi}_j$ and the approximating function $f(\mathbf{x})$ is represented as a sum of these RBFs, each associated with a different reference point $\bm{\xi}_j$, and weighted by a coefficient $c_j$ which has to be determined.

It can be seen that we get an overdetermined linear system of equations for the given dataset:
\begin{equation}
\label{eq:rbfAp}
\begin{aligned}
\displaystyle h_i &= f(\mathbf{x}_i)=\sum_{j=1}^M c_j \phi(\|\mathbf{x}_i-\bm{\xi}_j\|)\\
&=\sum_{j=1}^M c_j \phi_{i,j}\qquad i=1,\dots,N \textrm{.}
\end{aligned}
\end{equation}

\noindent The linear system of equations (\ref{eq:rbfAp}) can be represented in a matrix form as:
\begin{equation}
\label{eq:matE}
\mathbf{Ac=h}\textrm{,}
\end{equation}

\noindent where the number of rows is $N \gg M$ and $M$ is the number of unknown weights $[c_1,\dots, c_M]^T$, i.e. the number of reference points. Equation (\ref{eq:matE}) represents system of linear equations:

\begin{equation}
\label{eq:matLSE}
\begin{pmatrix}
  \phi_{1,1} & \cdots & \phi_{1,M} \\
  \vdots & \ddots & \vdots \\
  \phi_{i,1} & \cdots & \phi_{i,M} \\
  \vdots & \ddots & \vdots \\
  \phi_{N,1} & \cdots & \phi_{N,M}
 \end{pmatrix}
 \begin{pmatrix}
  c_1 \\
  \vdots\\
  c_M
 \end{pmatrix}=
  \begin{pmatrix}
  h_1 \\
  \vdots\\
  h_i\\
  \vdots\\
  h_N
 \end{pmatrix}\textrm{.}
\end{equation}
\\
\noindent The presented system is overdetermined, i.e. the number of equations $N$ is higher than the number of variables $M$. This linear system of equations can be solved by the least squares method as $\mathbf{A}^T\mathbf{A}\mathbf{c}=\mathbf{A}^T\mathbf{h}$ or singular value decomposition, etc.

\section{RBF Approximation for Large Data}\label{sec:big}
In practice, the real datasets contain a large number of points which results into high memory requirements for storing the matrix $\mathbf{A}$ of the overdetermined linear system of equations (\ref{eq:matLSE}). For example when we have dataset contains $3,000,000$ points, number of reference points is $10,000$ and double precision floating point is used then we need 223.5~GB memory for storing the matrix $\mathbf{A}$ of the overdetermined linear system of equations (\ref{eq:matLSE}). Unfortunately, we do not have an unlimited capacity of RAM memory and therefore calculation of unknown weights $c_j$ for RBF approximation would be prohibitively computationally expensive due to memory swapping, etc. In this section, a proposed solution to this problem is described.

In \cref{sec:RBF}, it was introduced that overdetermined system of equations can be solved by the least squares method. For this method the $M\times M$ square matrix:
\begin{equation}
\mathbf{B}=\mathbf{A}^T\mathbf{A}
\end{equation}
\noindent is to be determined. Advantages of matrix $\mathbf{B}$ are that it is a symmetric matrix and moreover only two vectors of length $N$ are needed to determine of one entry, i.e.:
\begin{equation}
b_{ij}=\sum_{k=1}^N \phi_{ki}\cdot \phi_{kj}\textrm{,}
\end{equation}
\noindent where $b_{ij}$ is the entry of the matrix $\mathbf{B}$ in the $i-$th row and $j-$th column. 

\begin{figure}[htb]
  \centering\includegraphics[width=.5\linewidth]{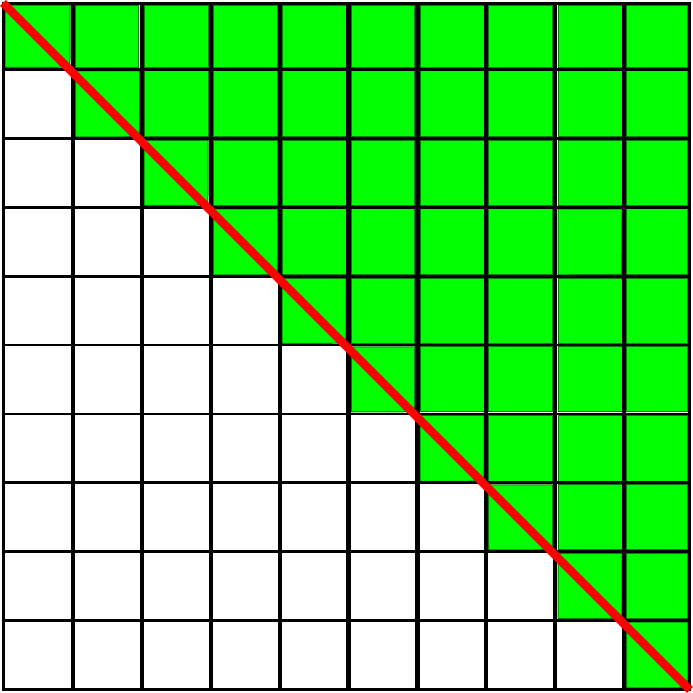}
  \caption{\label{fig:blockMatrix} $M\times M$ square matrix which is partitioned into $M_B\times M_B$ blocks. Main diagonal of matrix is represented by red color and illustrates the symmetry of matrix. Blocks, which must be computed, are represented by green color.}
\end{figure}

To save memory requirements and data bus (PCI) load block operations with matrices are used. Based on the above properties of the matrix $\mathbf{B}$, only the upper triangle of this matrix is computed. Moreover the matrix is partitioned into $M_B\times M_B$ blocks, see \cref{fig:blockMatrix}, and the calculation is performed sequentially for each block:
\begin{equation}
\begin{gathered}
\mathbf{B}_{kl}=(\mathbf{A}_{*,k})^T (\mathbf{A}_{*,l})\\[0.5em]
k = 1,\dots, \frac{M}{M_B}, \qquad l = k,\dots, \frac{M}{M_B} \textrm{,}
\end{gathered}
\end{equation}
where $\mathbf{B}_{kl}$ is sub-matrix in the $k-$th row and $l-$th column and $\mathbf{A}_{*,k}$ is defined as:
\begin{equation}
\displaystyle \mathbf{A}_{*,k}=\begin{pmatrix}
  \phi_{1,(k-1)\cdot M_B + 1} & \cdots & \phi_{1,k\cdot M_B} \\
  \vdots & \ddots & \vdots \\
  \phi_{i,(k-1)\cdot M_B + 1} & \cdots & \phi_{i,k\cdot M_B} \\
  \vdots & \ddots & \vdots \\
  \phi_{N,(k-1)\cdot M_B + 1}  & \cdots & \phi_{N,k\cdot M_B}
 \end{pmatrix}\textrm{.}
\end{equation}

The size of block $M_B$ is chosen so that $M_B$ is multiple of $M$ and there is no swapping, i.e.:
\begin{equation}
M_B \cdot(M_B + 2\cdot N)\cdot {prec}< \textrm{size of RAM [B],}
\end{equation}
where ${prec}$ is size of data type in bytes.

\section{Experimental results}
The presented modification of the RBF approximation method has been tested on synthetic and real data. Let us introduce results for two real datasets. 

The first dataset was obtained from LiDAR data of the Serpent Mound in Adams County, Ohio\footnotemark[1]. The second dataset is LiDAR data of the Mount Saint Helens in Skamania County, Washington\footnotemark[1]. Each point of these datasets is associated with its elevation. Summary of the dimensions of terrain for the given datasets is in \cref{tab:elevation}.
\footnotetext[1]{\url{http://www.liblas.org/samples/}}

\begin{table}[htb]
\caption{\label{tab:elevation} Summary of the dimensions of terrain for tested datasets. Note that one feet [ft] corresponds to 0.3048 meter [m].}
\renewcommand{\arraystretch}{1.2}
\begin{center}
\begin{tabular}{| c | r | r |}
\hline
\textbf{Dimensions} &  \multicolumn{1}{c|}{\textbf{Serpent Mound}} & \multicolumn{1}{c|}{\textbf{St. Helens}} \\ \hline \hline
\textbf{number of points} & $3,265,110$ & $6,743,176$\\ \hline
\textbf{lowest point [ft]} &  $166.7800$  & $3,191.5269$ \\ \hline
\textbf{highest point [ft]} & $215.4800$ & $8,330.2219$\\ \hline
\textbf{width [ft]}& $1,085.1199$ & $26,232.3696$\\ \hline
\textbf{length [ft]}& $2,698.9601$ & $35,992.6861$ \\ \hline

\end{tabular}
\end{center}
\end{table}

For experiments, two different radial basis functions have been used, see \cref{tab:Functions}. Shape parameters $\alpha$ for used RBFs were determined experimentally with regard to the quality of approximation and they are presented in \cref{tab:shapeParam}. Note that value of shape parameter $\alpha$ is inversely proportional to range of datasets. 

\begin{table}[htb]
\caption{\label{tab:Functions} Used RBFs }
\renewcommand{\arraystretch}{1.6}
\begin{center}
\begin{tabular}{|c | c | c | }
\hline
\textbf{RBF} & \textbf{type} & $\boldsymbol{\phi}(\bm{r})$ \\ \hline \hline
Gaussian RBF & global & $e^{-(\alpha r)^2}$ \\ \hline
Wendland's $\phi_{3,1}$ & local & $(1-\mathnormal{\alpha} r)_{+}^4(4\alpha r + 1)$ \\ \hline
\end{tabular}
\end{center}
\end{table}

\begin{table}[htb]
\caption{\label{tab:shapeParam} Experimentally determined shape parameters $\alpha$ for used RBFs}
\renewcommand{\arraystretch}{1.2}
\begin{center}
\begin{tabular}{| c | c | c |}
\hline
\multirow{2}{*}{\textbf{RBF}}  & \multicolumn{2}{c|}{\textbf{shape parameter}}  \\ \cline{2-3}
 &  \textbf{Serpent Mound} & {\textbf{St. Helens}} \\ \hline \hline
Gaussian RBF &  $\alpha = 0.05$  & $\alpha = 0.0004$ \\ \hline
Wendland's $\phi_{3,1}$ & $\alpha = 0.01$  & $\alpha = 0.0001$ \\ \hline
\end{tabular}
\end{center}
\end{table}

The set of reference points equals the subset of the given dataset for which we determine the RBF approximation. Moreover, the distribution of reference points is uniform and the set of reference points has a cardinality $10,000$ in both experiments.

\begin{figure*}
\centering\includegraphics[width=\textwidth]{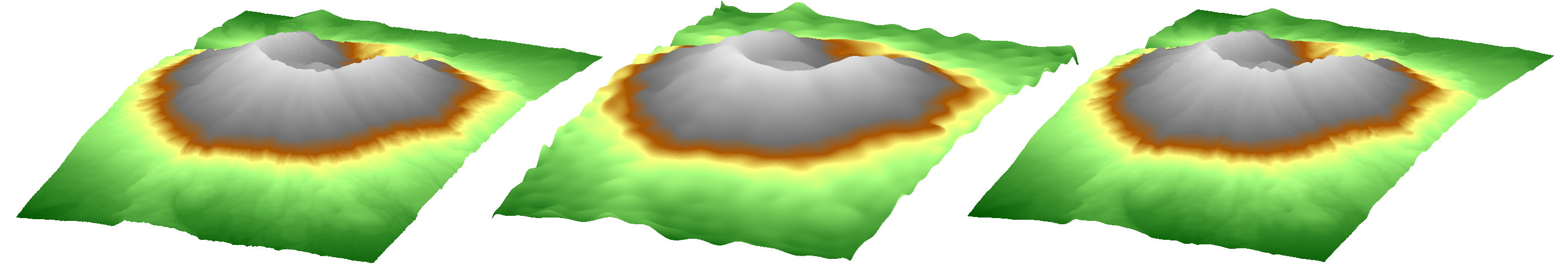}
\subfloat[\label{fig:Helens_orig}Original]{\includegraphics[width=0.3\textwidth]{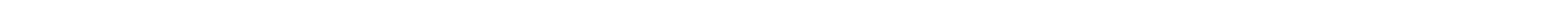}}
\subfloat[\label{fig:Helens_gauss}Gaussian RBF, $\alpha=0.0004$]{\includegraphics[width=0.3\textwidth]{img/blank.png}}
\subfloat[\label{fig:Helens_w31}Wendland's RBF $\phi_{3,1}$, $\alpha=0.0001$]{\includegraphics[width=0.3\textwidth]{img/blank.png}}\\
\centering\subfloat[\label{fig:Serpent_orig}Original]{\includegraphics[width=0.3\textwidth]{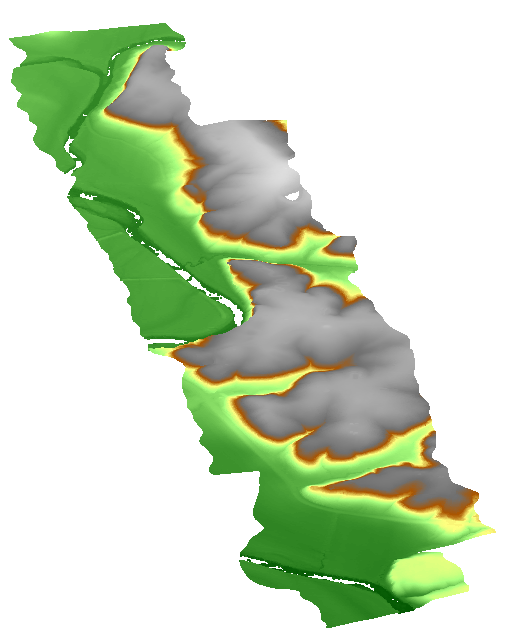}}
\subfloat[\label{fig:Serpent_gauss}Gaussian RBF, $\alpha=0.05$]{\includegraphics[width=0.3\textwidth]{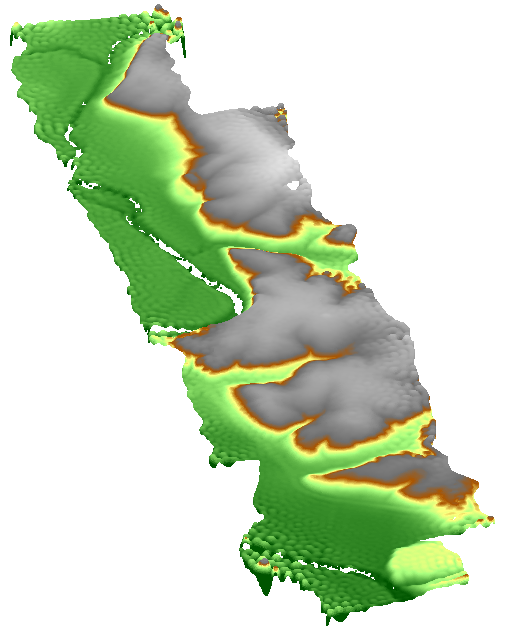}}
\subfloat[\label{fig:Serpent_w31}Wendland's RBF $\phi_{3,1}$, $\alpha=0.01$]{\includegraphics[width=0.3\textwidth]{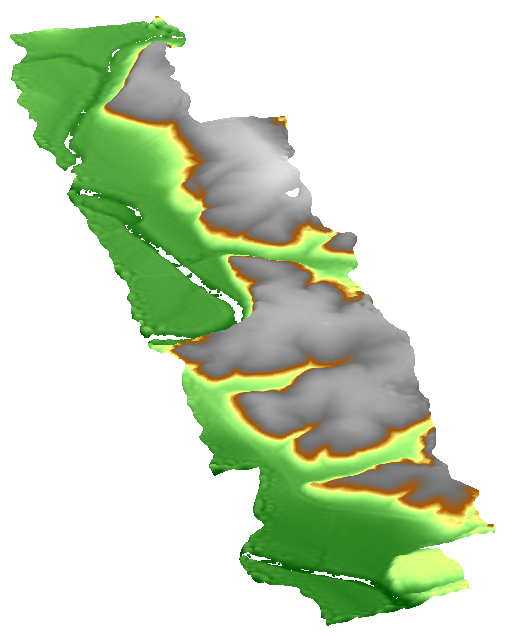}}
\vspace{1em}
\label{fig:HelensAndSerpent}\caption{Serpent Mound in Adams County, Ohio (top) and Mount Saint Helens is Skamania County, Washington (bottom)}
\end{figure*}

\begin{figure}[htb]
  \centering\includegraphics[width=\linewidth]{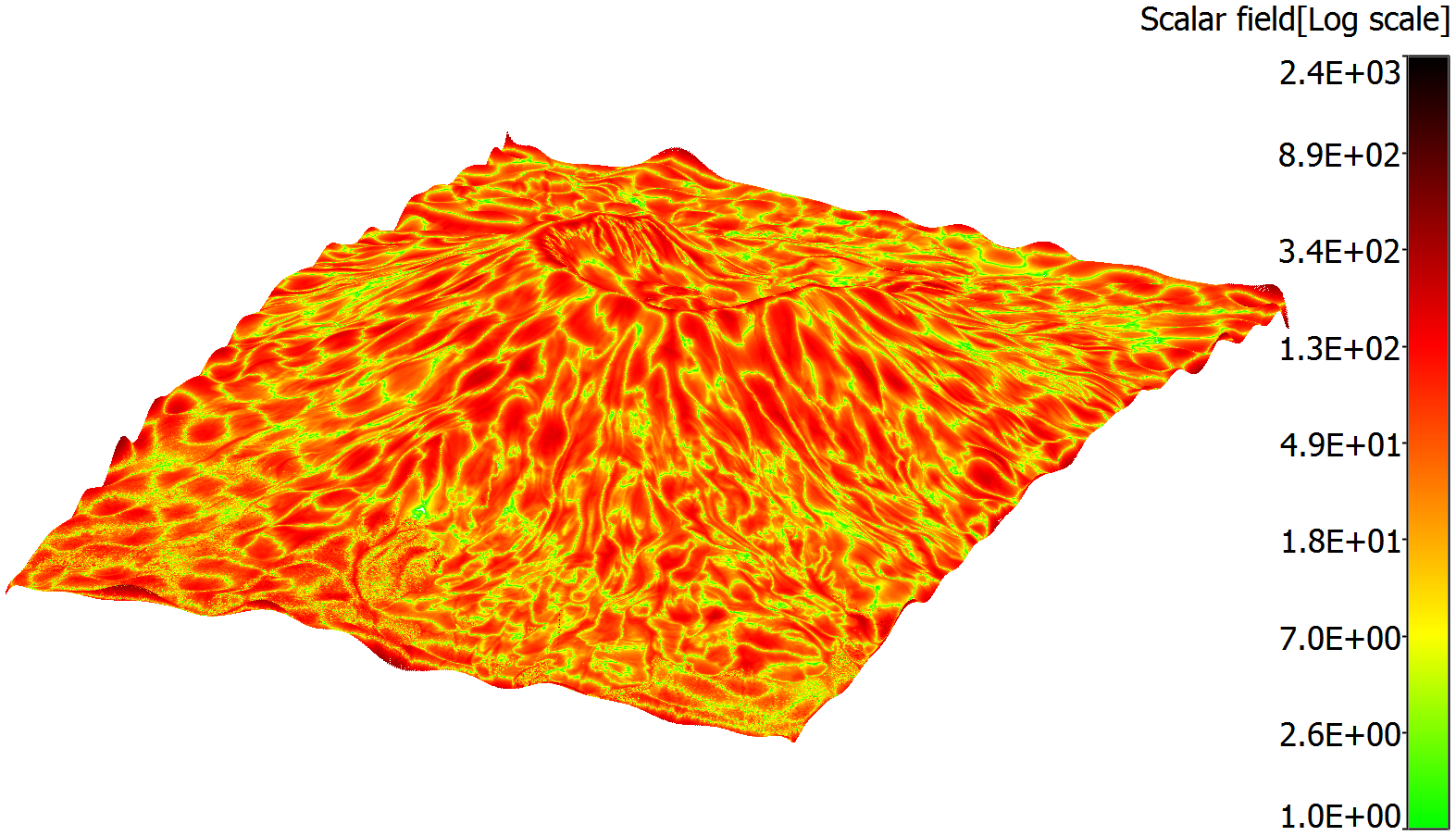}
  \caption{\label{fig:Helens_err_gaus} Approximation of Mount Saint Helens with $10,000$ global Gaussian basis functions with shape parameter $\alpha= \nobreak0.0004$ false-colored by magnitude of error.}
\end{figure}

\begin{figure}[htb]
  \centering\includegraphics[width=\linewidth]{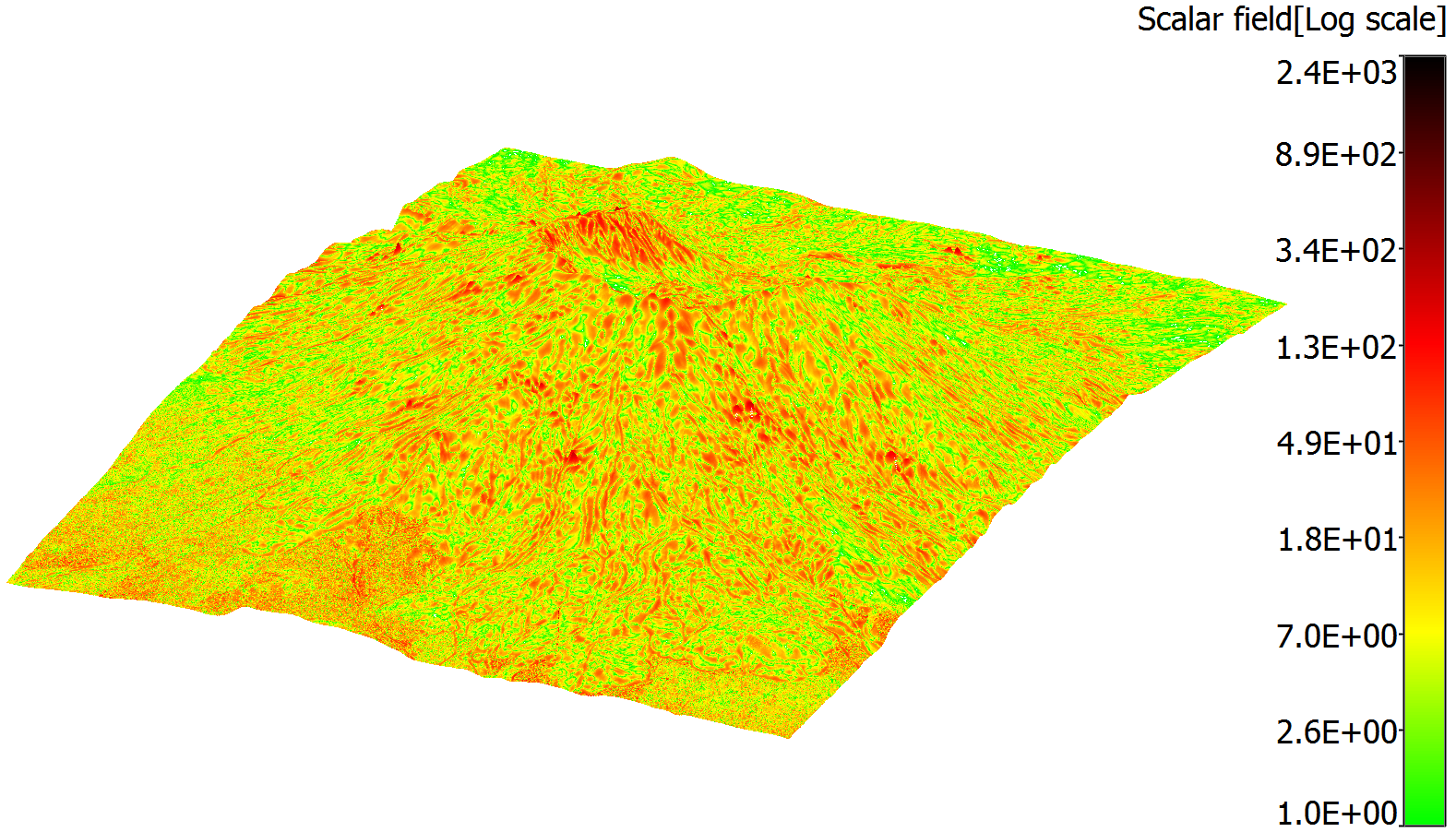}
  \caption{\label{fig:Helens_err_w31} Approximation of Mount Saint Helens with $10,000$ local Wendland's $\phi_{3,1}$ basis functions with shape parameter $\alpha=\nobreak0.0001$ false-colored by magnitude of error.}
\end{figure}

Approximation of Mount Saint Helens for both BRFs and its original are shown in \cref{fig:Helens_orig}-\ref{fig:Helens_w31}.
 In \cref{fig:Helens_gauss} can be seen that the RBF approximation with the global Gaussian RBFs cannot preserve the sharp rim of a crater. Further, visualization of magnitude of error at each point of the original points cloud is presented in \cref{fig:Helens_err_gaus} and \cref{fig:Helens_err_w31}. It can be seen that the RBF approximation with the global Gaussian RBFs returns worse result than RBF approximation with local Wendland's $\phi_{3,1}$ basis functions in terms of the error. In \cref{tab:err} can be seen the value of mean absolute error, its deviation and mean relative error for both approximations.

\begin{figure}[htb]
  \centering\includegraphics[width=.9\linewidth]{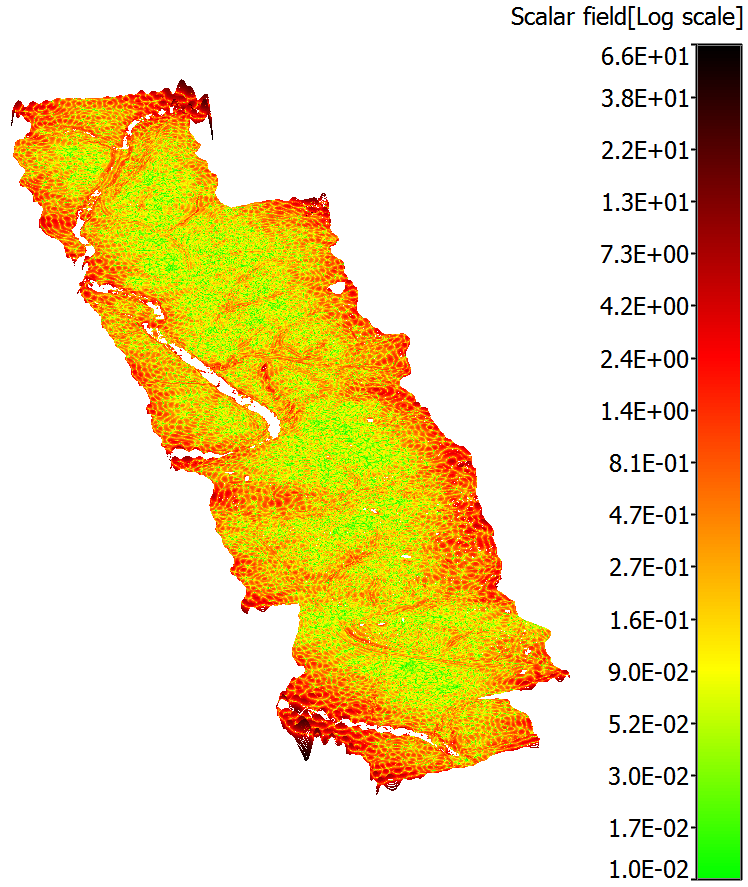}
  \caption{\label{fig:Serpent_err_gaus} Approximation of the Serpent Mound with $10,000$ global Gaussian basis functions with shape parameter $\alpha= \nobreak0.05$ false-colored by magnitude of error.}
\end{figure}

\begin{figure}[t!]
  \centering\includegraphics[width=.9\linewidth]{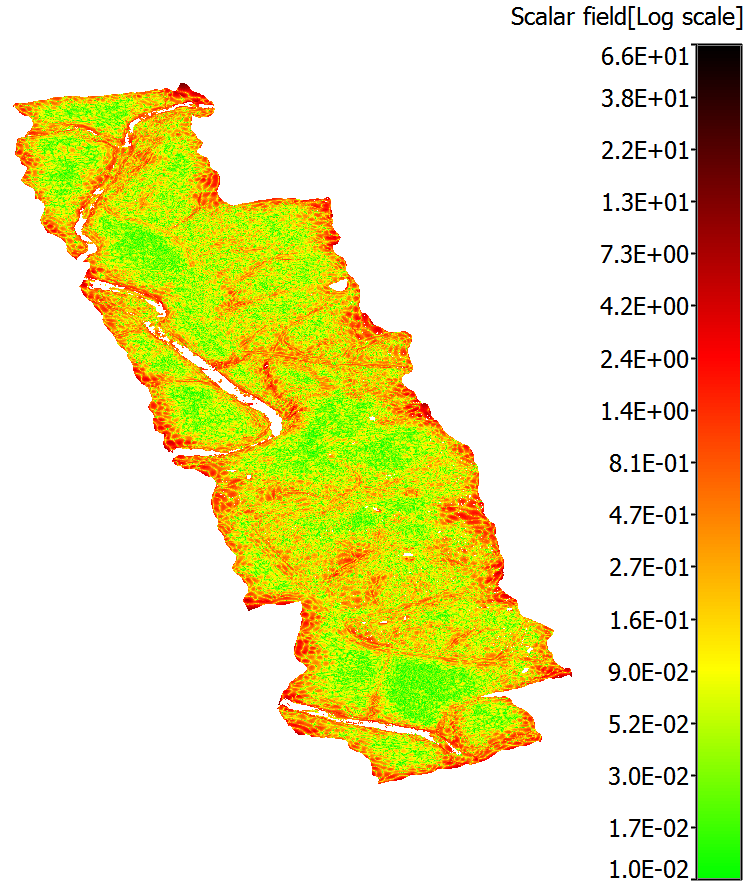}
  \caption{\label{fig:Serpent_err_w31} Approximation of the Serpent Mound with $10,000$ local Wendland's $\phi_{3,1}$ basis functions with shape parameter $\alpha=\nobreak0.01$ false-colored by magnitude of error.}
\end{figure}

Results of the RBF approximation for Serpent Mound and its original are shown in \cref{fig:Serpent_orig}-\ref{fig:Serpent_w31}. It can be seen that the approximation using local Wendland's $\phi_{3,1}$ basis function (\cref{fig:Serpent_w31}) returns again better result than approximation using the global Gaussian RBF (\cref{fig:Serpent_gauss}) in terms of the error. It is also seen in \cref{fig:Serpent_err_gaus} and \cref{fig:Serpent_err_w31} where magnitude of error at each point of original points cloud is visualized.
Moreover, we can see that the highest errors occur on the boundary of terrain, which is a general problem of RBF methods. Value of mean absolute error, its deviation and mean relative error due to elevation for both used RBFs are again mentioned in \cref{tab:err}.

Mutual comparison both datasets in terms of the mean relative error (\cref{tab:err}) indicates that mean relative error for Serpent Mount is smaller than for Mount Saint Helens. It is caused by the presence of vegetation, namely forest, in LiDAR data of the Mount Saint Helens. This vegetation operates in our RBF approximation as noise and therefore the resulting mean relative error is higher.

The implementation of the RBF approximation has been performed in Matlab and tested on PC with the following configuration:
\begin{itemize}
\item CPU: Intel\textregistered~Core\texttrademark~i7-4770 (4$\times$ 3.40GHz + hyper-threading),\\
\item memory: 32 GB RAM,\\
\item operating system Microsoft Windows 7 64bits.
\end{itemize}

\begin{figure}[b!]
\centering
\includegraphics[width=0.91\linewidth]{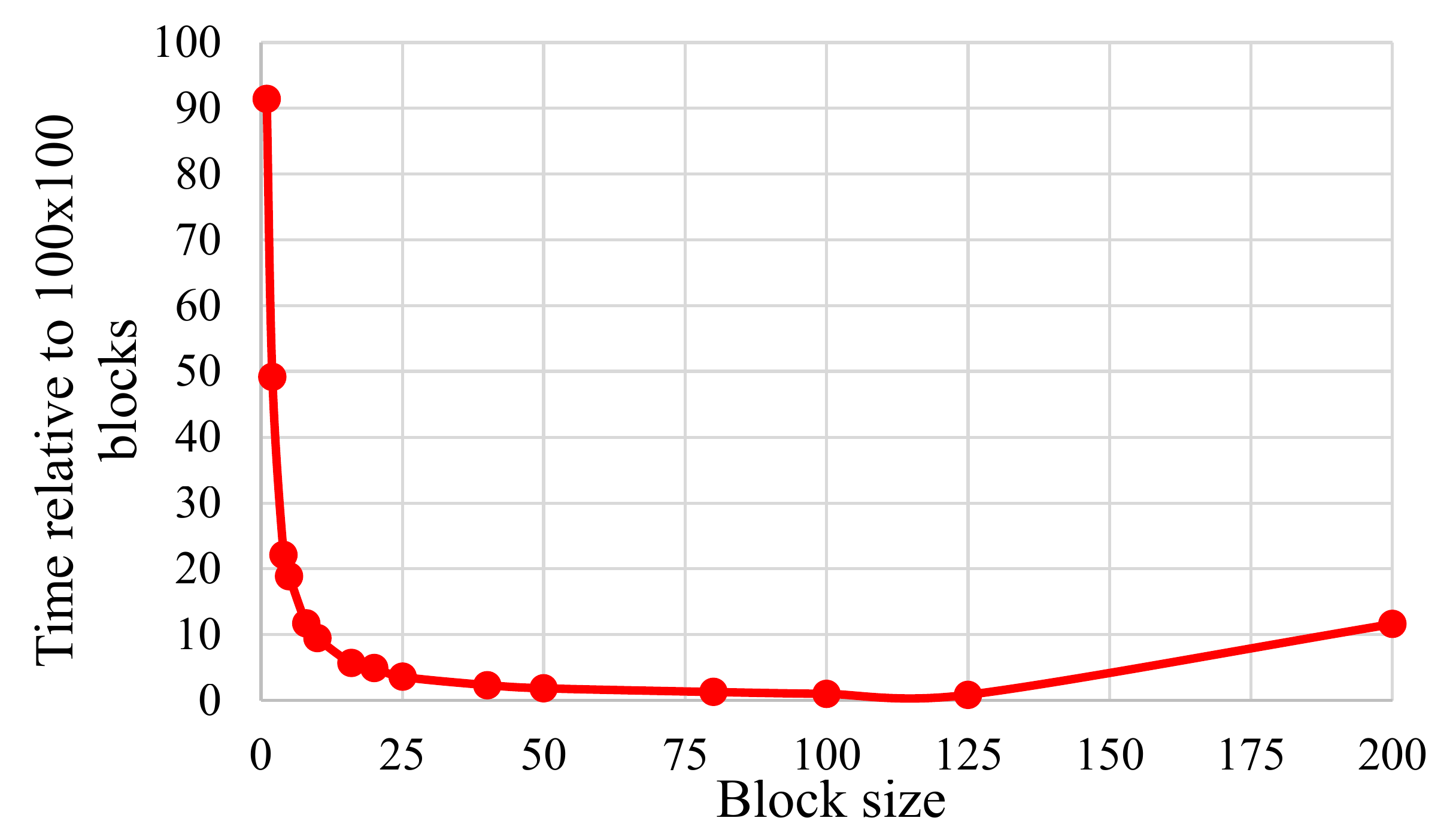}\\
\includegraphics[width=0.91\linewidth]{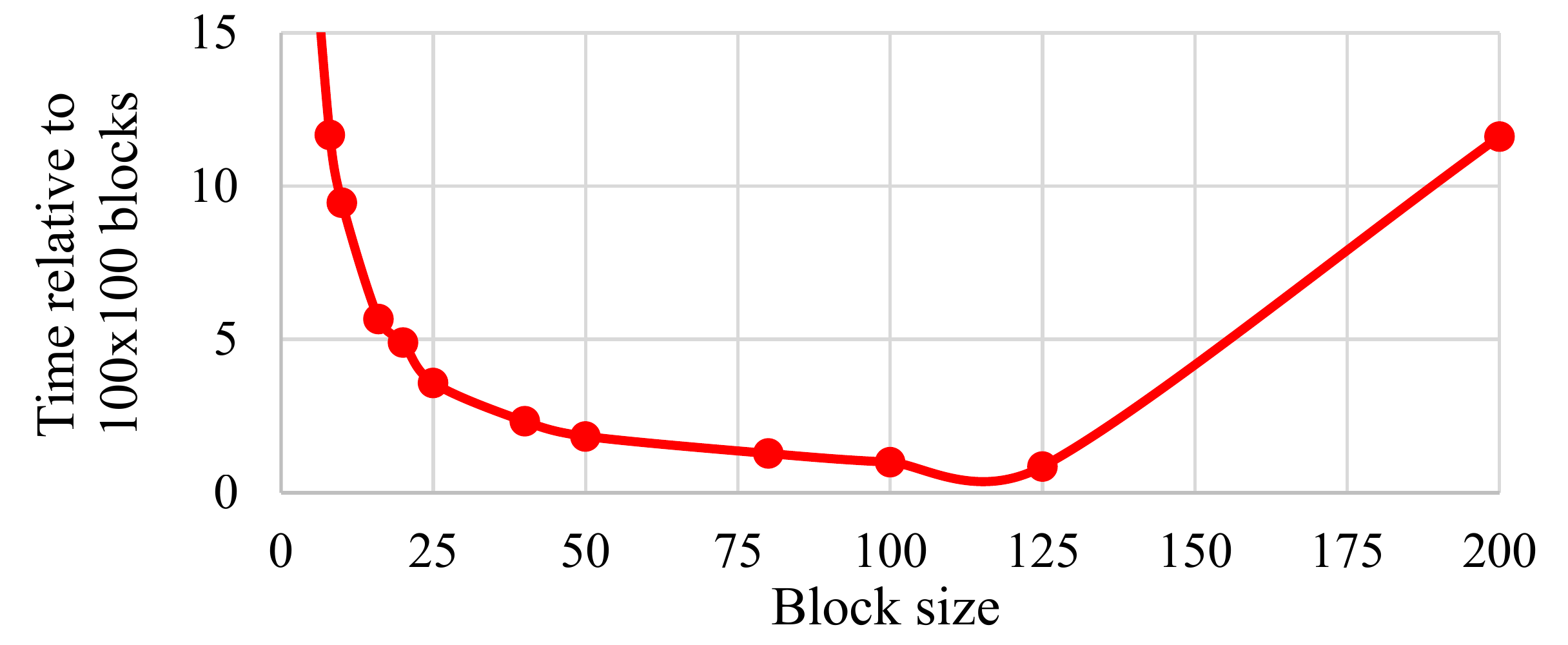}\\[-1em]
\caption{\label{fig:time} Time performance for approximation of the Serpent Mound depending on the block size. The times are presented relative to the time for $100\times 100$ blocks.}
\end{figure}

 \noindent For the approximation of the Serpent Mound with $10,000$ local Wendland's $\phi_{3,1}$ basis function with shape parameter $\alpha=0.01$ the running times  for different sizes of blocks were measured. These times were converted relative to the time for $100\times 100$ blocks and are presented in \cref{fig:time}. We can see that for the approximation matrix which is partitioned into small blocks (i.e. smaller than $25\times25$ blocks) the time performance is large. This is caused by overhead costs. On the other hand, for the approximation matrix which is partitioned into large blocks (i.e. larger than $125\times125$ blocks) the running time begins to grow above the permissible limit due to memory swapping. 

\begin{table*}
\centering
\caption{The RBF approximation error for testing datasets and different radial basis functions. Note that one feet [ft] corresponds to 0.3048 meter [m].}
\label{tab:err}
\renewcommand{\arraystretch}{1.2}
\begin{tabular}{ |l |  r | r | r |r|}
\hline
\multicolumn{1}{|c|}{\multirow{2}{*}{\textbf{Error}}}& \multicolumn{2}{c|}{\textbf{Serpent Mound}} & \multicolumn{2}{c|}{\textbf{St. Helens}}\\ \cline{2-5}
\multicolumn{1}{|c|}{}& \multicolumn{1}{c|}{\textbf{Gaussian RBF}} & \multicolumn{1}{c|}{\textbf{Wendland's} $\phi_{3,1}$} & \multicolumn{1}{c|}{\textbf{Gaussian RBF}} & \multicolumn{1}{c|}{\textbf{Wendland's} $\phi_{3,1}$} \\ \hline\hline
\textbf{mean absolute error [ft]} & 0.4477 & 0.2289 & 44.4956& 12.1834\\
\textbf{deviation of error [ft]}      & 1.4670 & 0.1943 & 680.3659& 169.2800\\
\textbf{mean relative error [\%]} & 0.0024 & 0.0012 & 0.0087& 0.0023\\ \hline
\end{tabular}
\end{table*}

\section{Conclusions}
This paper presents a new approach to the RBF approximation of large datasets. The proposed approach uses symmetry of matrix and partitioning matrix into blocks, thus preventing memory swapping. The experiments made proved that the proposed approach is able to determine the RBF approximation for large dataset. Moreover, from the experimental results we can see that use of a local RBFs is better than global RBFs, if data are sufficiently sampled. Futher, it is obvious that approximation using the global Gaussian RBFs has problems with the preservation of sharp edges. The experiments made also proved that RBF methods have problems with the accuracy of calculation on the boundary of an object, which is a well known property, and the magnitude of the RBF approximation error is influenced by the presence of a noise.

For the future work, the RBF approximation method can be explored in terms of lower sensitivity to noise, more accurate calculation on the boundary or better approximation of sharp edges and improvements of the computational cost without loss of approximation accuracy.

\section*{Acknowledgments}
The authors would like to thank their colleagues at the University of West Bohemia, Plzen, for their discussions and suggestions, and also anonymous reviewers for the valuable comments and suggestions they provided. The research was supported by MSMT CR projects LH12181 and SGS~2016-013.


\bibliographystyle{eg-alpha}

\bibliography{egbibsample}

\newcommand{\etalchar}[1]{$^{#1}$}
\begin{thebibliography}{\uppercase{CBC{\etalchar{*}}01}}

\bibitem[CBC{\etalchar{*}}01]{carr2001}
\textsc{Carr J.~C., Beatson R.~K., Cherrie J.~B., Mitchell T.~J., Fright W.~R.,
  McCallum B.~C., Evans T.~R.}:
\newblock Reconstruction and representation of 3d objects with radial basis
  functions.
\newblock In \emph{Proceedings of the 28th Annual Conference on Computer
  Graphics and Interactive Techniques, {SIGGRAPH} 2001, Los Angeles,
  California, USA, August 12-17, 2001} (2001), pp.~67--76.

\bibitem[Dar00]{Darve2000195}
\textsc{Darve E.}:
\newblock The fast multipole method: Numerical implementation.
\newblock \emph{Journal of Computational Physics 160}, 1 (2000), 195--240.

\bibitem[Fas07]{fasshauer2007}
\textsc{Fasshauer G.~E.}:
\newblock \emph{{M}eshfree {A}pproximation {M}ethods with {MATLAB}}, vol.~6.
\newblock World Scientific Publishing Co., Inc., River Edge, NJ, USA, 2007.

\bibitem[Har71]{Hardy1971}
\textsc{Hardy R.~L.}:
\newblock {M}ultiquadratic {E}quations of {T}opography and {O}ther {I}rregular
  {S}urfaces.
\newblock \emph{Journal of Geophysical Research 76} (1971), 1905--1915.

\bibitem[HSfY15]{Hon2015}
\textsc{Hon Y.-C., Sarler B., fang Yun D.}:
\newblock Local radial basis function collocation method for solving
  thermo-driven fluid-flow problems with free surface.
\newblock \emph{Engineering Analysis with Boundary Elements 57} (2015), 2 -- 8.
\newblock \{RBF\} Collocation Methods.

\bibitem[LCC13]{Li2013}
\textsc{Li M., Chen W., Chen C.}:
\newblock The localized \{RBFs\} collocation methods for solving high
  dimensional \{PDEs\}.
\newblock \emph{Engineering Analysis with Boundary Elements 37}, 10 (2013),
  1300 -- 1304.

\bibitem[PRF14]{pepper2014}
\textsc{Pepper D.~W., Rasmussen C., Fyda D.}:
\newblock A meshless method using global radial basis functions for creating
  3-d wind fields from sparse meteorological data.
\newblock \emph{Computer Assisted Methods in Engineering and Science 21}, 3-4
  (2014), 233--243.

\bibitem[PS11]{pan2011two}
\textsc{Pan R., Skala V.}:
\newblock A two-level approach to implicit surface modeling with compactly
  supported radial basis functions.
\newblock \emph{Eng. Comput. (Lond.) 27}, 3 (2011), 299--307.

\bibitem[Ska13]{skala2013}
\textsc{Skala V.}:
\newblock {F}ast {I}nterpolation and {A}pproximation of {S}cattered
  {M}ultidimensional and {D}ynamic {D}ata {U}sing {R}adial {B}asis {F}unctions.
\newblock \emph{WSEAS Transactions on Mathematics 12}, 5 (2013), 501--511.

\bibitem[Ska15]{Skala2015}
\textsc{Skala V.}:
\newblock Meshless interpolations for computer graphics, visualization and
  games.
\newblock In \emph{Eurographics 2015 - Tutorials, Zurich, Switzerland, May 4-8,
  2015} (2015), Zwicker M., Soler C., (Eds.), Eurographics Association.

\bibitem[SPN13]{nedved2013}
\textsc{Skala V., Pan R., Nedved O.}:
\newblock Simple 3d surface reconstruction using flatbed scanner and 3d print.
\newblock In \emph{{SIGGRAPH} Asia 2013, Hong Kong, China, November 19-22,
  2013, Poster Proceedings} (2013), {ACM}, p.~7.

\bibitem[SPN14]{skala2014}
\textsc{Skala V., Pan R., Nedved O.}:
\newblock Making 3d replicas using a flatbed scanner and a 3d printer.
\newblock In \emph{Computational Science and Its Applications - {ICCSA} 2014 -
  14th International Conference, Guimar{\~{a}}es, Portugal, June 30 - July 3,
  2014, Proceedings, Part {VI}} (2014), vol.~8584 of \emph{Lecture Notes in
  Computer Science}, Springer, pp.~76--86.

\bibitem[TO02]{Turk2002}
\textsc{Turk G., O'Brien J.~F.}:
\newblock Modelling with implicit surfaces that interpolate.
\newblock \emph{{ACM} Trans. Graph. 21}, 4 (2002), 855--873.

\bibitem[Wen06]{wendland2006}
\textsc{Wendland H.}:
\newblock Computational aspects of radial basis function approximation.
\newblock \emph{Studies in Computational Mathematics 12} (2006), 231--256.

\end{thebibliography}

\newpage

\end{document}